\documentclass[%
 reprint,
 amsmath,amssymb,
 aps,
]{revtex4-2}

\usepackage{graphicx}
\usepackage{dcolumn}
\usepackage{bm}
\usepackage{hyperref}
\hypersetup{colorlinks=blue,linkcolor=blue,citecolor=blue}

\usepackage{float}
\begin{document}

\preprint{APS/123-QED}

\title{Anticipating Tipping Points for Disordered Traffic: Critical Slowing Down on the Onset of Congestion}

\author{Shankha Narayan Chattopadhyay}
\email{shankha.21maz0008@iitrpr.ac.in}
\author{Arvind Kumar Gupta}%
 \email{akgupta@iitrpr.ac.in}
\affiliation{%
 Department of Mathematics, Indian Institute of Technology Ropar,
  Rupnagar, 140 001, Punjab, India 
}%

\date{\today}
\begin{abstract}
Regime shifts are frequently observed in intricate systems ranging from natural systems such as population collapse, ecosystems, marine ice instability, etc to man-made systems like vehicular traffic flow, disease outbreaks, etc. A number of statistical indicators, commonly referred to as early warning signals (EWS), have been proposed in order to predict these sudden transformations beforehand. These signals have received much acclaim for their general characteristics. In this study, we proposed a lattice hydrodynamic area occupancy model with passing to examine the dynamics of a heterogeneous disordered traffic system, which is particularly significant in the context of emerging economies where road infrastructure often lacks lane discipline. The creation of kink or chaotic jam has been analyzed theoretically and also demonstrated using numerical simulations. Additionally, we have established a suitable framework to address the forthcoming issue of traffic congestion by employing generic statistical measures. We have successfully proved through the use of simulated data that our measures exhibit sensitivity to the phase transition occurring between free flow and congested flow.
\end{abstract}

\maketitle

\section{Introduction}
In the last few decades due to globalization and rapid urbanization, the number of automobiles has been continuously growing, as a direct consequence traffic congestion has emerged as a momentous issue across the world. For this reason, the problem of traffic jams has attracted considerable attention from researchers, engineers, and policymakers \cite{jin2015understanding,tadaki2013phase} because of its high impact on human lifestyle as well as on the environment \cite{chin1996containing}. Developed countries have better infrastructure, but they are still struggling to tackle the traffic-related consequences. In developing countries, traffic management and control becomes an extremely difficult task, because not only the infrastructure is poorly designed but also the traffic is exceptionally heterogeneous and disordered.
Sophisticated technologies based on intelligent transport systems are helpful in implementing smart city traffic control management \cite{wang2023transportation}. The advancement of high-level computing efficiencies in recent times has prompted the use of data-assisted and data-driven methodologies in order to enhance the performance of multimodal transport systems. These techniques aim to tackle the issues associated with traffic congestion analysis and forecasts \cite{wang2005real,huang2022adaptive}. Though these data-driven techniques are highly advantageous, their applicability remains a question in places where sufficient data availability is an issue \cite{kumar2015short}.  
The modeling of traffic flow is a key tool to examine the complex physical mechanism of congestion and to understand the interactions of vehicles on the roadway.
The pursuit to understand different complex transport phenomena, both natural as well as man-made, leads to the study of non-equilibrium statistical mechanics \cite{schadschneider2010stochastic}.
The first known mathematical model regarding traffic flow was perhaps proposed by Lighthill and Whitham \cite{lighthill1955kinematic} and Richards \cite{richards1956shock} independently.
Currently, there has been extensive research conducted on traffic issues at several levels. These investigations encompass:  microscopic models such as car-following models 
\cite{ma2021nonlinear}, cellular automata models 
\cite{li2001cellular}
and macroscopic models such as continuum models 
\cite{gupta2005analyses}.
Macroscopic continuum models consider traffic flow to be analogous to fluid flow and are comprised of coupled partial differential equations 
\cite{zhai2022continuum,zhang2021extended}
Higher order non-equilibrium continuum models are useful to design suitable control strategies, as they deal with the fundamental variables density, velocity, and flow, however, they suffer from several drawbacks, for example, high cost of computation, numerical instability, failure to express simple real-life traffic situations like passing, individual driver's characteristics, etc.\\
Experimental evidence suggests that traffic congestion emerges from free flow when some of the drivers vary their speeds and the distance between the consecutive vehicles changes, vehicles begin to bunch up and a cluster comes to a near standstill situation, stop and go traffic waves appear and the jam spreads backward over time 
\cite{sugiyama2008traffic}.
At the hydrodynamic level, utilizing a lattice framework, the appearance of nonlinear stop-and-go traffic waves has been explored through numerical simulations \cite{nagatani1998modified}. This lattice hydrodynamic (LH) modeling approach describes the occurrence of different traffic waves in terms of kink-antikink and soliton density waves; like a microscopic model, and analyzes the collective dynamics of traffic; like a macroscopic model. Recently, these LH models got attention from researchers, due to their advantages, simplicity, flexibility, and high computational efficacy 
\cite{nagatani1999chaotic,redhu2015jamming}.\\
From the point of view of physics, traffic jams are the result of a dynamical phase transition due to sudden changes in the parameters influencing the overall traffic system. A nearly free flow can suddenly meet congestion if the mean car density of a segment surpasses some critical threshold because of some external influences like accidents, rough driving, bad weather conditions, etc.  Abrupt and often irreversible transitions can be observed not only in traffic flow framework \cite{ghadami2021stability}, but also in a wide variety of systems including social networks \cite{bhandary2023rising}, and cellular biological processes \cite{sarkar2019anticipating} to name a few.
These so-called `critical transitions' (CT) are typically brought on by a gradual change in external conditions that quietly bring the system into the vicinity of a `tipping point' (TP). TPs are actually threshold values at which the system dynamics switch to some alternate state due to some stochastic perturbations. As a stable system approaches its state of linear stability, the real component of the dominating eigenvalue associated with the linearized system converges toward zero. Consequently, a local bifurcation is approached, and the first-order approximation of the potential function gets flattened, as a consequence the system becomes more vulnerable to perturbations \cite{lenton2013environmental,lenton2011early}. The presence of perturbations results in extended periods of transient recovery towards the equilibrium state in close proximity to the approaching TP. Despite the seemingly unpredictable nature of an approaching TP, there are certain mathematical features that a system exhibits on this journey toward it. Those features can provide us with clues as to the risk of an imminent critical transition. As a dynamical system approaches a TP (i.e. traffic congestion in our study), the system dynamics as a whole exhibit a gradual decrease in speed, which leads to the phenomenon, known as `critical slowing down' (CSD). CSD is at the heart of many of the `early warning signals' (EWSs) developed today \cite{boers2021critical}. These EWSs are useful to anticipate these critical transitions before they actually arrive. Monitoring simple statistical measures of the time series may be selected as indicators of approaching instability \cite{dai2012generic,carpenter2006rising,dakos2012methods}. These statistical measures tell us how close the system is to instability by monitoring the system behavior in the free flow regime. A sharp hike in the indicators is a symptom of upcoming instability. However, there exist questions about the statistical robustness of EWSs \cite{ditlevsen2010tipping} and further studies have revealed that there are critical transitions with no prior EWSs \cite{ditlevsen2010tipping,hastings2010regime}. We can use a combination of two or more statistical properties of the data to get robust early warning indicator of approaching bifurcation\cite{ditlevsen2010tipping}.\\For a lane-based traffic system without overtaking, using a higher order homogeneous macroscopic continuum model \cite{jiang2002new}, a recent study has demonstrated that the examination of the generic EWSs can effectively serve as a reliable indicator for predicting the occurrence of upstream stop-and-go traffic congestion \cite{ghadami2021stability}. But, in real-life scenarios passing or overtaking is a vitally important part of highway traffic. Now question arises about the robustness of generic EWSs in the presence of overtaking. In a developing country like India, the majority of roads consist a single lane on which vehicles of different types of static and dynamic characteristics interact. When there is availability of lateral space small vehicles respond more quickly than large vehicles due to their high sensitivity, and accelerate-decelerate rapidly to utilize the available space for overtaking. It is indeed interesting to analyze how this type of heterogeneity affects the overall disordered traffic dynamics. Our objective of this study is to analyze the following concerns: (I) How does the heterogeneity and passing affect the appearance of congestion for a unidirectional disordered traffic flow? (II) Are generic EWSs sensitive to the rapid phase transition between free flow and congested flow on a highway, consisting of heterogeneous and non-lane-based disordered traffic? (III) Are EWSs applicable for anticipating all types of regime shifts that happen due to changes in heterogeneity?
\begin{figure*}[!t] 
    \centering
       \includegraphics[width=17cm,height=5.5cm]{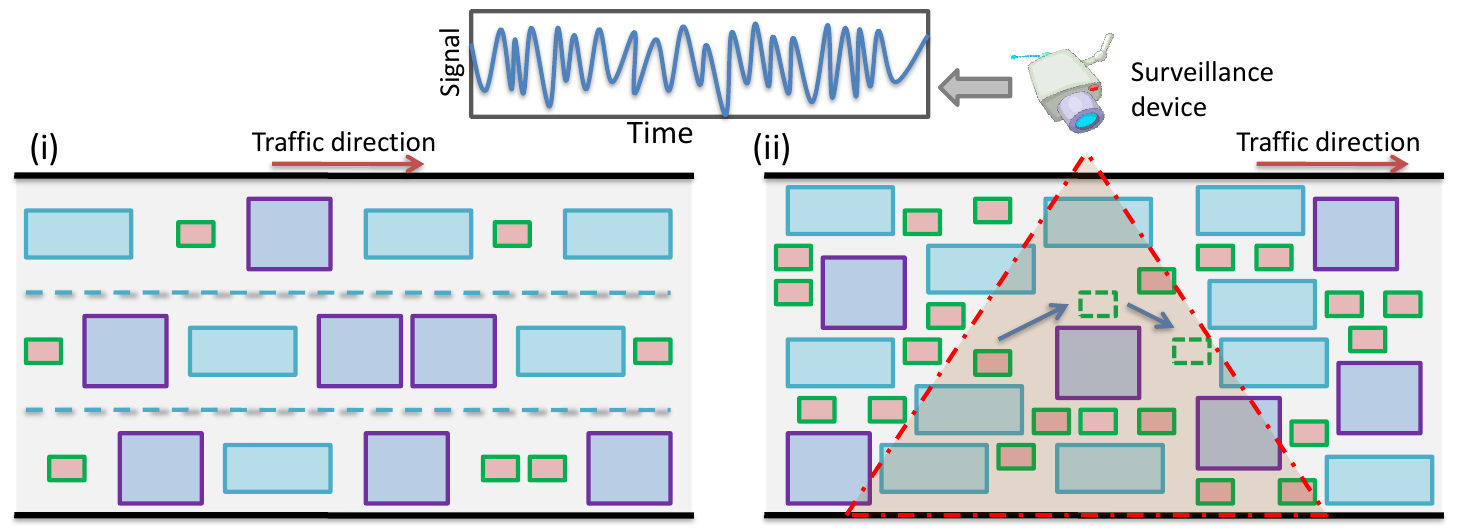}
    \\ (a)\\
        \includegraphics[width=8.5cm,height=4.5cm]{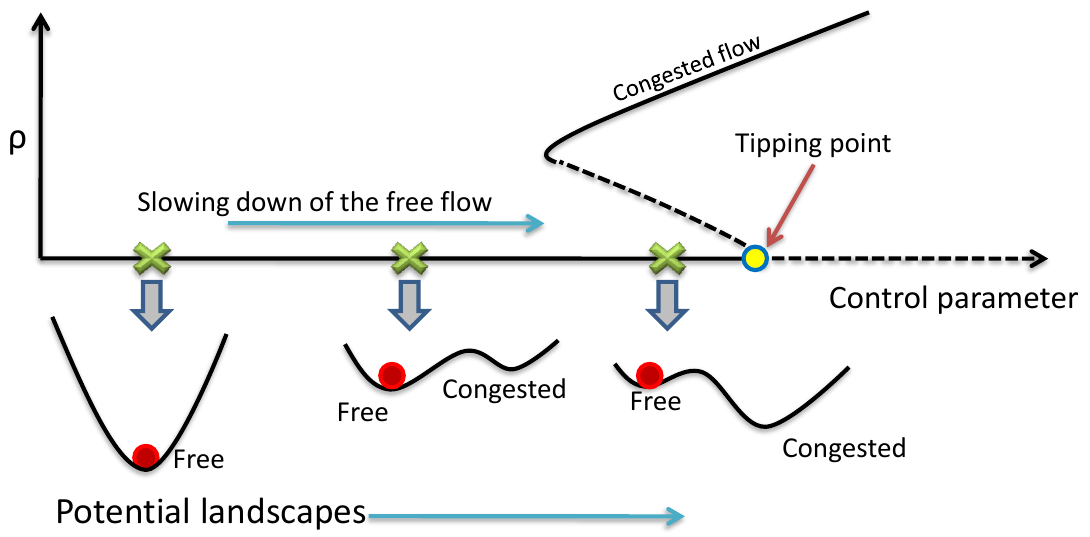}\quad
    \includegraphics[width=8.5cm,height=4.5cm]{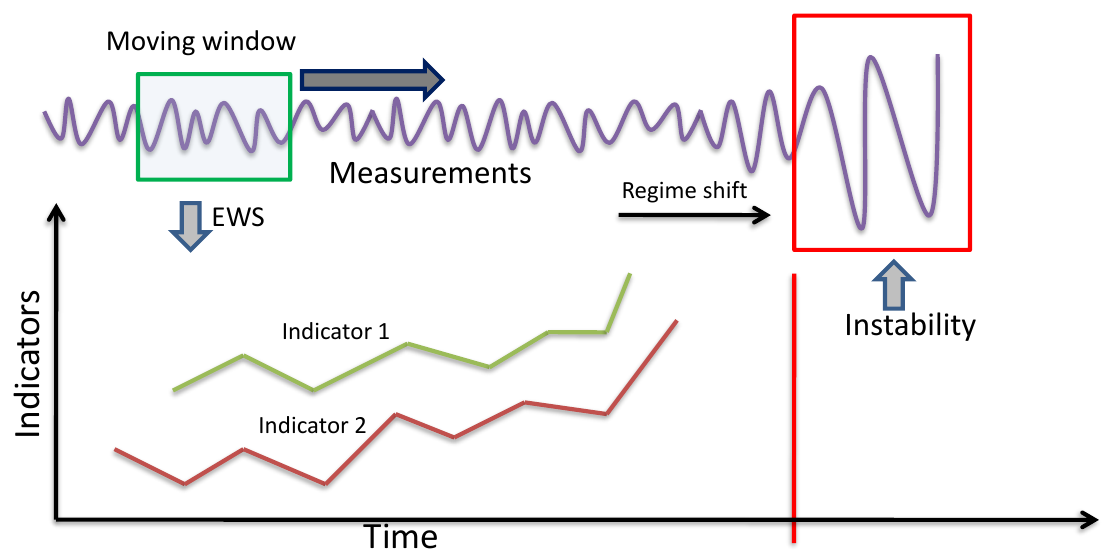}\\
 (b) \hspace{8cm} (c)\\
\caption{(a)There are two common types of heterogeneous traffic systems found on roads, which involve several vehicles with distinct physical properties and behavioral characteristics. (i) The concept of an ideal lane-based traffic system involves the implementation of strict lane discipline, where vehicles of various types adhere to designated lanes. (ii) On the other hand, a realistic traffic system may not always maintain perfect lane order, allowing vehicles of different sizes to travel alongside each other and enabling overtaking maneuvers by utilizing the lateral width of the road. In the context of a disordered traffic system, a surveillance device is employed to capture and quantify the movement of traffic, afterward converting it into a signal that can be subjected to analysis. (b) This graphic presents a systematic representation of the critical slowing down phenomenon observed in a traffic flow system when it approaches a tipping point. As the tipping point is near, the possible landscape becomes increasingly level. The stable attractor is represented by the local minima of the potential well, while the state of the system is shown by the presence of the red ball.  As the tipping point approaches, the ball becomes increasingly susceptible to external disturbances. The acquisition of moving window statistical indicators is derived from the analysis of observed traffic patterns. An abrupt increase in the indications before the occurrence of the regime transition serves as an early indication of the imminent regime shift. 
}
  \label{fig1} 
\end{figure*}
\section{Methodology} \label{Metho}
\subsection{Traffic Flow Modelling}
 To extract data from a heterogeneous disordered traffic system, a non-lane based area occupancy LH model allowing passing, has been employed. The model comprises two types of differential-difference equations, namely the continuity equation and the evolution equation. We consider $m$ types of vehicles moving on a closed loop road having width $W$ and the road is divided into $L$ distinct lattice sites. Let, $\rho_{l,j}(t)$ and $v_{l,j}(t)$, respectively are the density and velocity for $l$-th ($1 \leq l \leq m$) type of vehicle, at $j$-th ($1\leq j \leq L$) lattice site at time $t$. In the absence of any on/off ramp, the total number of vehicles will remain conserved and as a consequence,  the continuity equation for $l$-th type of vehicle in lattice structure can be expressed in the following manner
\begin{equation}
    \partial_t\bigl(\rho_{l,j}(t)\bigr)+ \rho_0\bigl(q_{l,j}(t)-q_{l,j-1}(t)\bigr) = 0,\\
    \label{eq:1}
 \end{equation}
where   $q_{l,j}(t)=\rho_{l,j}(t)v_{l,j}(t)$ is the corresponding flow of $l$-th type of vehicle on $j$-th site at time $t$ and $\rho_0$ is the average density of the entire road segment. \\
To ensure the coherence of the system, it is necessary to provide an additional flow equation that incorporates two fundamental characteristics of traffic flow: the tendency of vehicles to achieve the equilibrium flow, and the passing/overtaking effect. Passing occurs from the $j$-th site to $j+1$-th site when the traffic current at $j$-th site exceeds the current at $j+1$-th site, and the extent of passing is directly proportional to the difference between the optimal currents at $j$-th site and $j+1$-th site. When a vehicle attempts to overtake another, it assesses the necessary space based on its own shape and size before making a choice, and therefore passing rate can be taken as a constant $\gamma_l$ for $l$-th  type of vehicles. In a typical
traffic situation, vehicles of varying sizes have the ability to move alongside one another using the lateral space on a shared roadway (Fig. \ref{fig1}(a)(ii)). Hence, in order to accommodate the lack of such lane discipline, along with the length it is necessary to take into account the entire width of the road section \cite{mohan2017heterogeneous}. In the context of heterogeneous traffic, where multiple types of vehicles coexist on the road, the metric used to quantify the extent of area occupancy is defined as $B=\sum_{l=1}^m\frac{c_lA_l}{W}$ where $c_l$ $\bigl(\sum_{l=1}^m c_l =1\bigr)$ and $A_l$  are respectively the fraction and area occupied by the $l$-th  type vehicle \cite{kaur2022analyses}. The velocity of a certain class of vehicle is influenced by the interaction with all types of vehicles, which is quantified by the total density $\rho_{j}(t)=\sum_{l=1}^m \rho_{l,j}(t)$ of all classes of vehicles at the $j$-th lattice site at time $t$. The integration of the impact resulting from the interaction between various types of vehicles, along with the influence of area occupancy, has been accomplished by incorporating the modified optimal velocity function  $V\bigl(\rho_j^*(t)\bigr)$ where $\rho^*_{j}(t)=B\rho_{j}(t)$. Now, the corresponding equation of flow for $l$-th type of vehicle with the additional effect of area occupancy allowing passing is given by
\begin{equation}\label{eq:2}
 \begin{split}
\partial_t\bigl(q_{l,j}(t)\bigr)&=a_l\Bigl[c_l\rho_0\frac{v_l^{max}}{2}V\bigl(\rho^*_{j+1}(t)\bigr)-q_{l,j}(t)\Bigr]\\
   &+a_l\gamma_l\rho_0\frac{v_l^{max}}{2}\Bigl[V\bigl(\rho^*_{j+1}(t)\bigr)-V\bigl(\rho^*_{j+2}(t)\bigr)\Bigr],
\end{split}
\end{equation}
 where  $a_l$ is the sensitivity and $v_l^{max}$ is the maximum allowed speed of $l$-th type of  vehicle. The modified version of the optimal velocity function \cite{kaur2022analyses} considering the additional effect of area occupancy is defined as
 \begin{equation}
V\bigl(\rho^*_j(t)\bigr)=\biggl[tanh\Bigl(\frac{2}{\rho_0}-\frac{\rho^*_{j}(t)}{\rho_0^2}-\frac{1}{\rho_c}\Bigr)+tanh\Bigl(\frac{1}{\rho_c}\Bigr)\biggr].
\label{eq:4}
\end{equation}
The function $V\bigl(\rho_j^*(t)\bigr)$ has the inflection point at $\rho_j^*(t)=\rho_c$ when $\rho_0=\rho_c$. In order to eliminate $v_{l,j}(t)$ from Eqs. (\ref{eq:1}) and (\ref{eq:2}), we discretize the partial derivatives $\partial_t\bigl(\rho_{l,j}(t)\bigr)$ and $\partial_t\bigl(q_{l,j}(t)\bigr)$ by first order forward difference scheme  with step size $\tau_l=\frac{1}{a_l}$. Next, the elimination of $v_{l,j}(t)$ in the system of Eqs. (\ref{eq:1}) and (\ref{eq:2}) through a one-to-one correspondence, for $l$-th type of vehicles we get the single density evolution equation as
\begin{equation}\label{eq:2v2}
\Delta_t\rho_{l,j}(t+\tau_l)=c_l\tau_l\rho_0^2 \frac{v_l^{max}}{2} \Bigl(\gamma_l\Delta^{2}V\bigl(\rho_j^*(t)\bigr)-\Delta V\bigl(\rho_j^*(t)\bigr)\Bigr),
    \end{equation}
    where
     $$\begin{cases} 
      \Delta_{t} \rho_{l,j}(t+\tau_l)=\rho_{l,j}(t+2\tau_l)-\rho_{l,j}(t+\tau_l),\\
       \Delta V\bigl(\rho_j^*(t)\bigr)=V\bigl(\rho_{j+1}^*(t)\bigr)-V\bigl(\rho_j^*(t)\bigr),\\
     \Delta^{2}V(\rho_j^*(t))=V\bigl(\rho_{j+2}^*(t)\bigr)-2V\bigl(\rho_{j+1}^*(t)\bigr)+V\bigl(\rho_{j}^*(t)\bigr).
   \end{cases}$$
To get the overall density equation for the mixed traffic, we add  Eq. (\ref{eq:2v2}) for $l=1$ to $l=m$,
 \begin{equation}\label{eq:2v3}
\sum_{l=1}^m \Delta_t\rho_{l,j}(t+\tau_l)=\rho_0^2M\Delta^{2}V\bigl(\rho_j^*(t)\bigr)-\rho_0^2K\Delta V\bigl(\rho_j^*(t)\bigr),
\end{equation}
  where   
     $K=\sum_{l=1}^m\frac{c_l\tau_lv_l^{max}}{2}$, and $M=\sum_{l=1}^m\frac{c_l\tau_l\gamma_l v_l^{max}}{2}.$
The primary objective of mixed traffic is to provide a comprehensive representation of traffic patterns on a global scale, rather than emphasizing specific attributes of individual vehicles. Hence, we introduce $\tau$ as the overall delay parameter such as $\tau_l=k_l\tau$, with $k_l \leq 1$ and $k_l < k_{l+1}$ as time delay associated with smaller cars is shorter compared to that of larger cars. Thus, the LHS of Eq. (\ref{eq:2v3}) can be rewritten as
 \begin{equation}\label{eq:2v4}
 \begin{split}
\sum_{l=1}^m \Delta_t \rho_{l,j}(t+\tau_l)&= \sum_{l=1}^m \rho_{l,j}\bigl((t+2\tau)+2\tau(k_l-1)\bigr) \\&- \sum_{l=1}^m \rho_{l,j}\bigl((t+\tau)+\tau(k_l-1)\bigr)
\end{split}
\end{equation}
Utilizing the Taylor series expansion for RHS of Eq. (\ref{eq:2v4})  and further neglecting terms involving derivatives in the expansion, from Eq. (\ref{eq:2v3}) we get the density evolution equation for $m$ types of vehicles in mixed traffic as follows
 \begin{equation}\label{eq:5v1}
\Delta_{t} \rho_j(t+\tau)=\tau\rho_0^2 \Bigl(D\Delta^{2}V\bigl(\rho_j^*(t)\bigr)-C\Delta V\bigl(\rho_j^*(t)\bigr)\Bigr),
 \end{equation}

where $$ \begin{cases} \Delta_{t} \rho_j(t+\tau)=\rho_j(t+2\tau)-\rho_j(t+\tau),\\
     C=\sum_{l=1}^m\frac{c_lk_lv_l^{max}}{2}, D=\sum_{l=1}^m\frac{c_lk_l\gamma_lv_l^{max}}{2}.\end{cases}$$
   
We define $\gamma=\frac{D}{C}$ and the mathematical model for heterogeneous disordered traffic without on/off ramp is given by
    \begin{equation}\label{eq:5v2}
\Delta_{t} \rho_j(t+\tau)=C\tau\rho_0^2 \Bigl(\gamma\Delta^{2}V\bigl(\rho_j^*(t)\bigr)-\Delta V\bigl(\rho_j^*(t)\bigr)\Bigr).
 \end{equation}
 When one type of vehicle moves in homogeneous traffic with passing but without the effect of area occupancy Eq. (\ref{eq:5v2}) yields \cite{nagatani1999chaotic} and therefore $\gamma=\frac{D}{C}$ serves like overall passing efficacy of the disordered mixed traffic. The preservation of the continuity equation is observed even in the scenario of heterogeneity and passing.\\
Further to tackle the additional effect of on/off ramp, we have considered the generation rate $g_j(t)$ as the difference of in and out flow at $j$-th lattice. $g_j(t)$ is zero when there is no ramp on the road otherwise at the location of an in/out ramp $g_j(t)=q_{in}\delta(j-j_{in})-q_{out}\delta(j-j_{out})$, where $\delta$ is the Kronecker delta function. While we have incorporated the essential characteristics of traffic flow on a unidirectional non-lane based heterogeneous disordered system in Eq. (\ref{eq:5v2}), it is necessary to consider the influence of drivers' errors in estimating road density and adjusting their speed based on the optimal velocity function. To accurately represent this phenomenon, the inclusion of spatiotemporal stochasticity is required. Spatiotemporal stochasticity can be represented by additive  Gaussian white noise $\xi_j(t)$ in space and time.  Finally, the stochastic difference equation representing heterogeneous disordered traffic with an on/off ramp is as follows
   \begin{equation}\label{eq:5}
   \begin{split}
\Delta_{t} \rho_j(t+\tau)&=C\tau\rho_0^2 \Bigl(\gamma\Delta^{2}V\bigl(\rho_j^*(t)\bigr)-\Delta V\bigl(\rho_j^*(t)\bigr)\Bigr)\\&+g_j(t)+\xi_j(t).
\end{split}
 \end{equation}
  \begin{figure}[!t] \centering
\includegraphics[width=4.25cm,height=4cm]{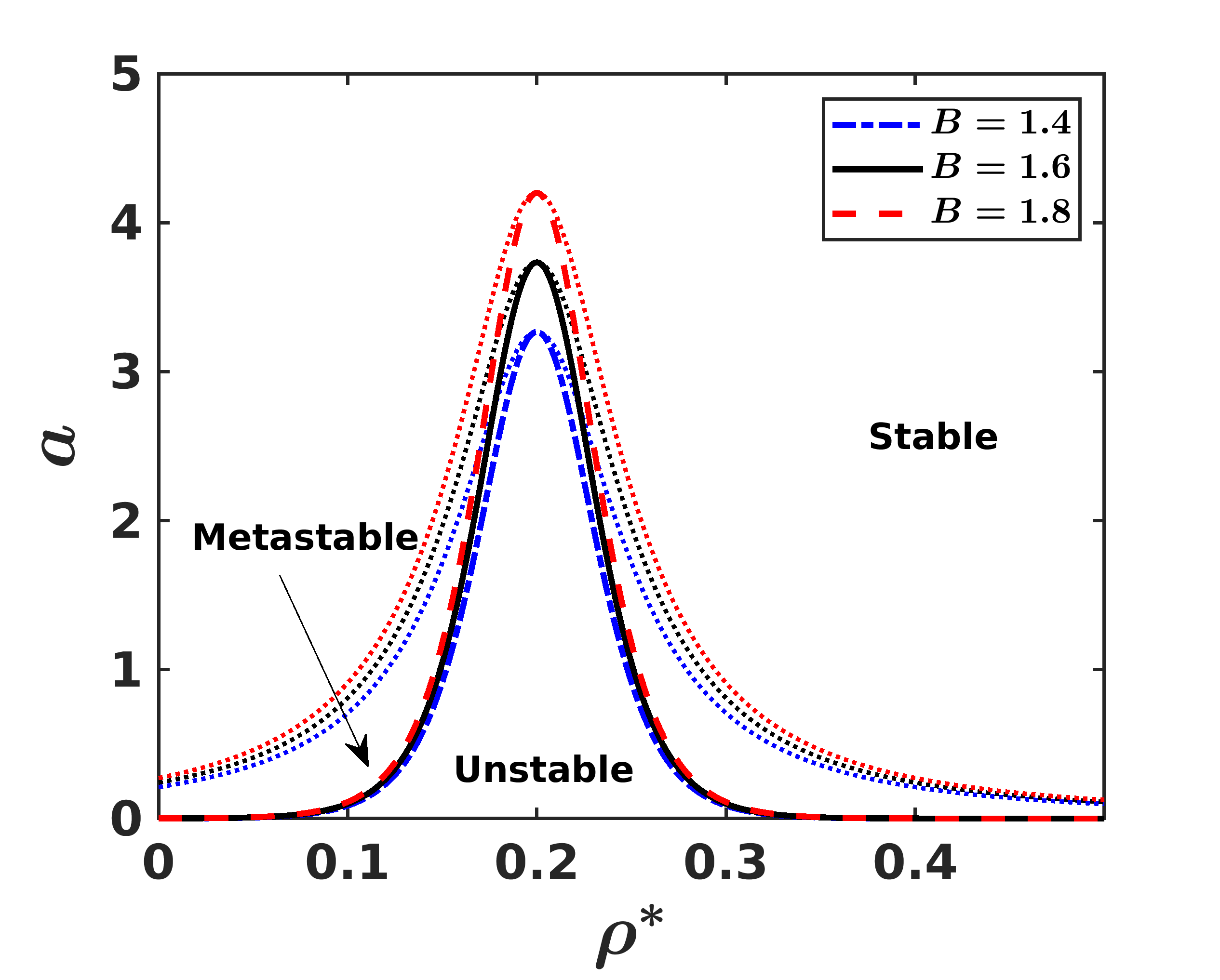}\hfill
   \includegraphics[width=4.25cm,height=4cm]{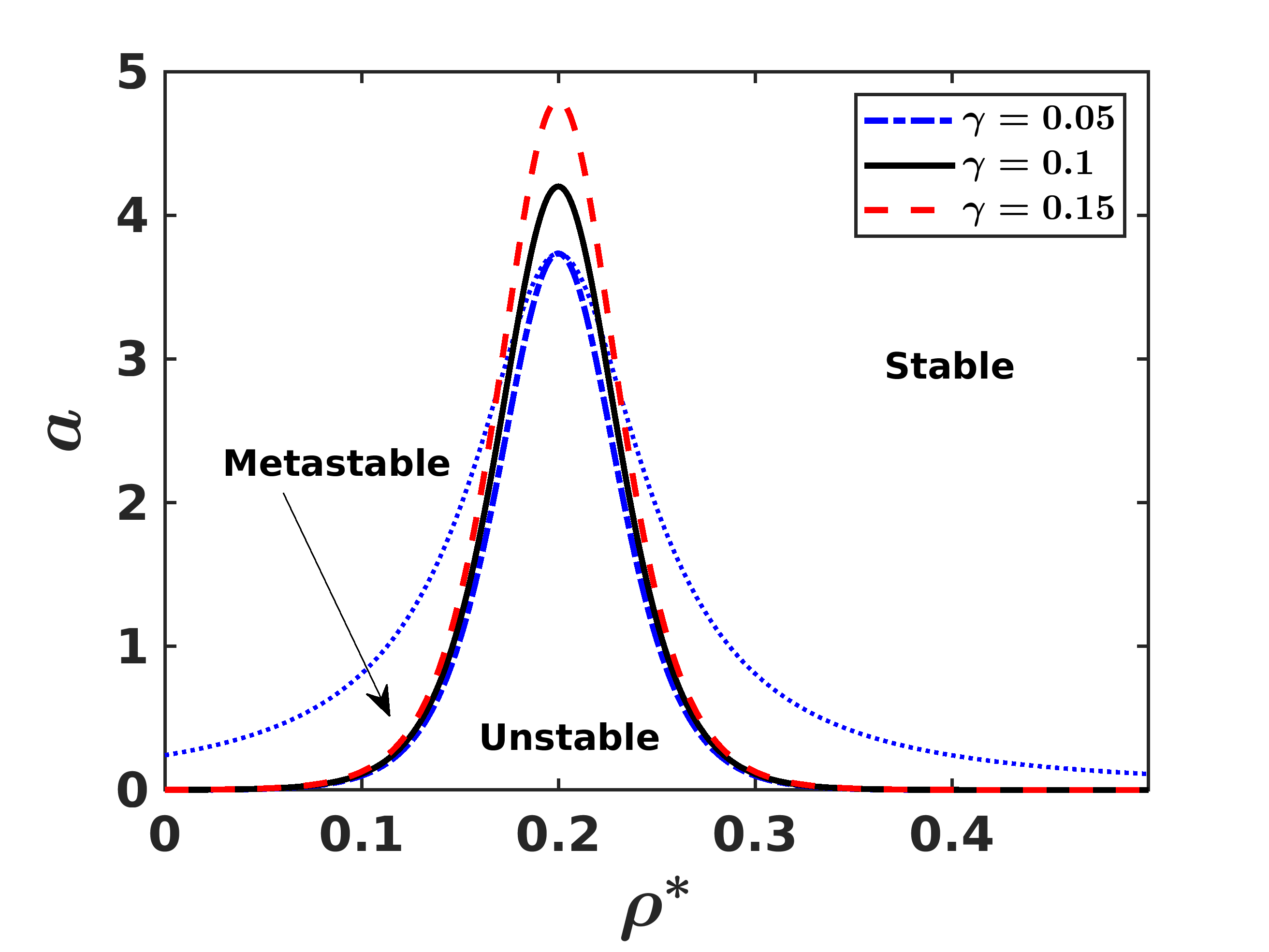}\\
  (a) \hspace{5cm} (b) \\
        \includegraphics[width=4.25cm,height=4cm]{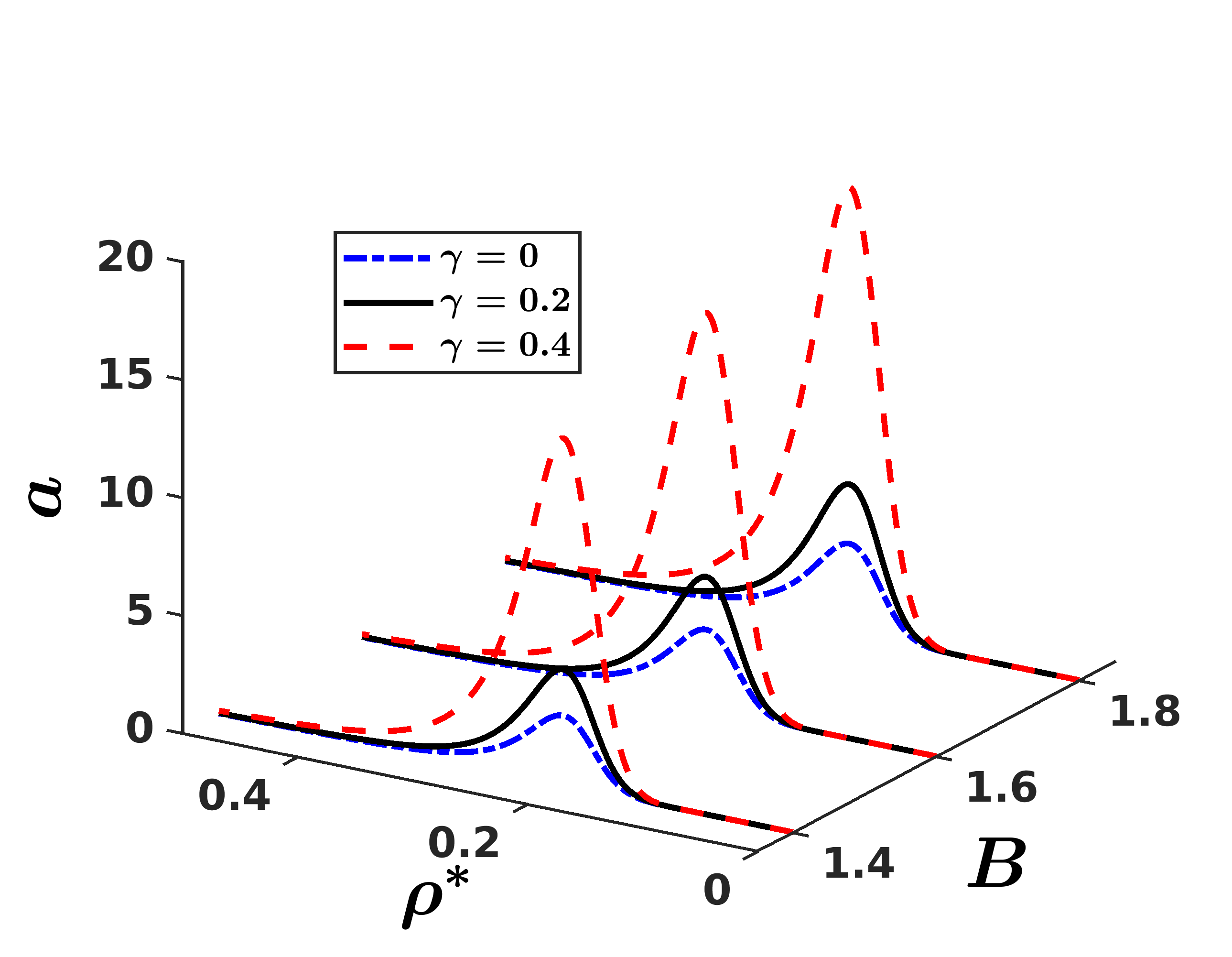}\hfill
   \includegraphics[width=4.25cm,height=4cm]{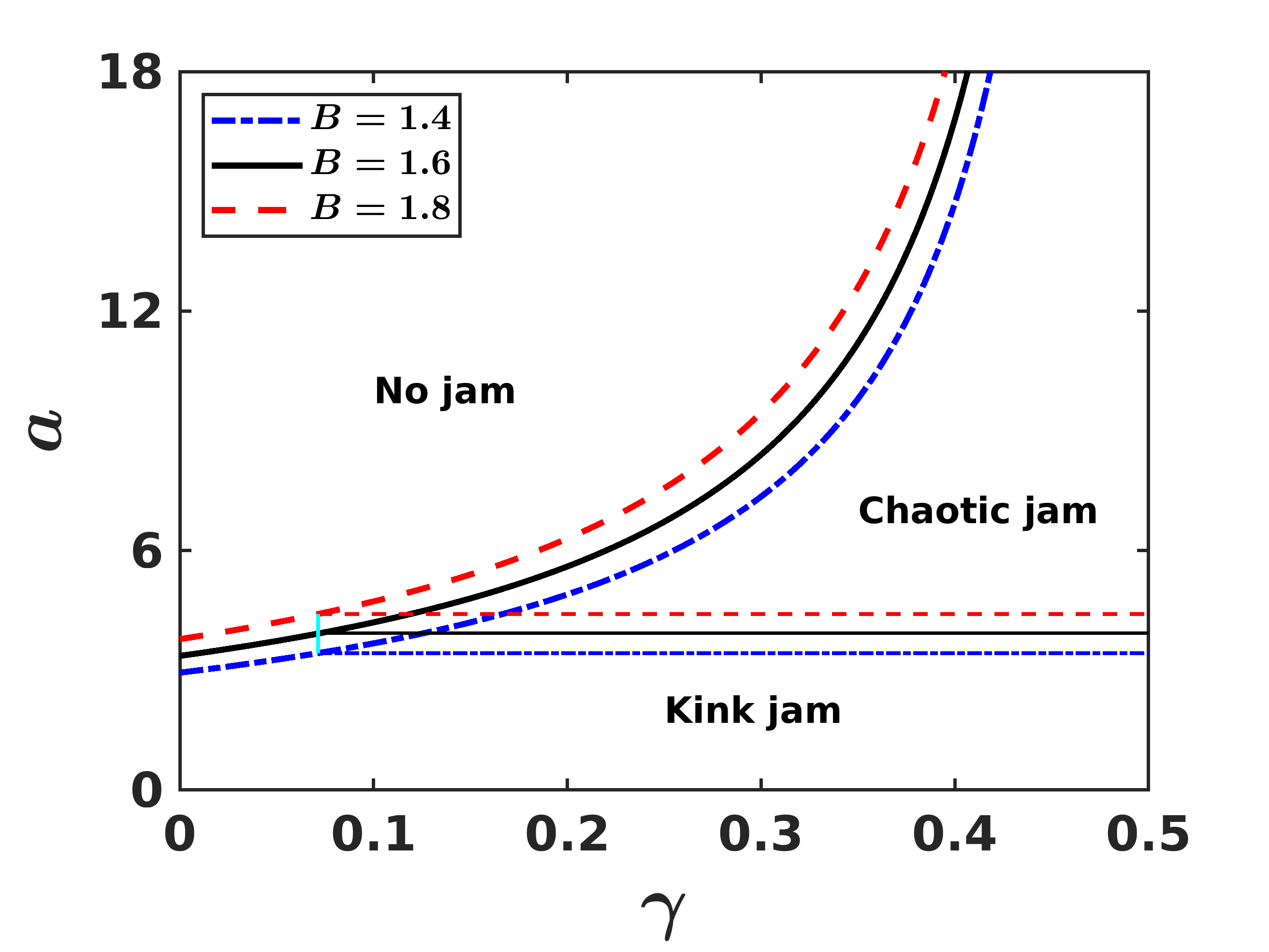}\\
    (c) \hspace{5cm} (d) \\
  \caption{ Phase diagrams of density $\rho^*_j(t)$$\bigl(=B\rho_j(t)\bigr)$ versus sensitivity $a$ for (a) different values of area occupancy factor $B$ when $C=0.7$ and $\gamma=0.05$ and (b) for different values of passing rate $\gamma$ when  $C=0.7$ and $B=1.6$. (c) Sensitivity diagram in parameter space $(\rho^*,B,a)$ space when $C=0.7.$ (d) Phase diagrams of passing rate $\gamma$ versus sensitivity $a$ for different values of area occupancy factor $B$ when $C=0.7$ and $\rho_0=0.2$.}
\label{fig2} 
\end{figure}
 \begin{figure*}[!t] 
    \centering
         \includegraphics[width=5.25cm,height=4.5cm]{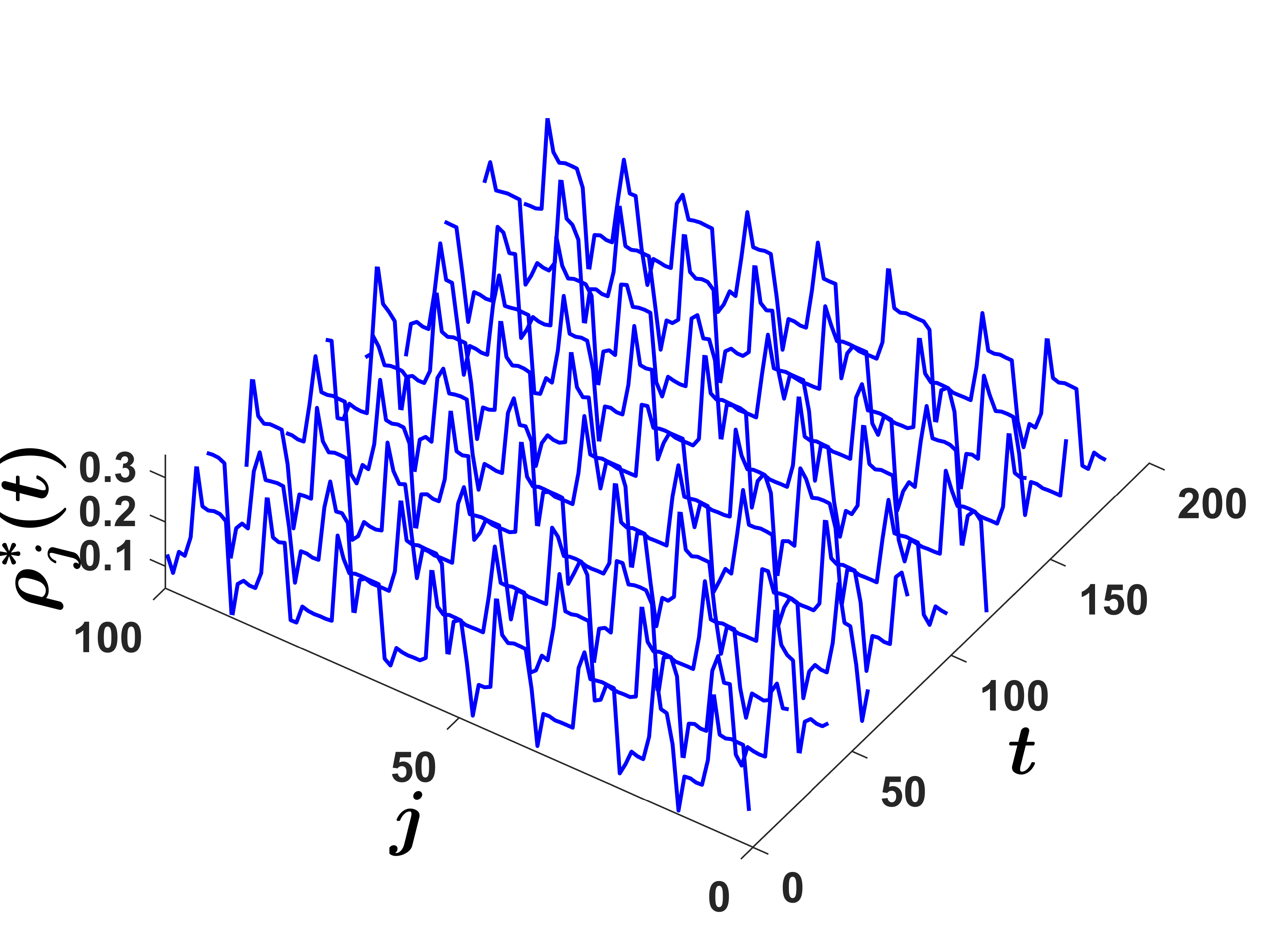}\quad
      \includegraphics[width=5.25cm,height=4.5cm]{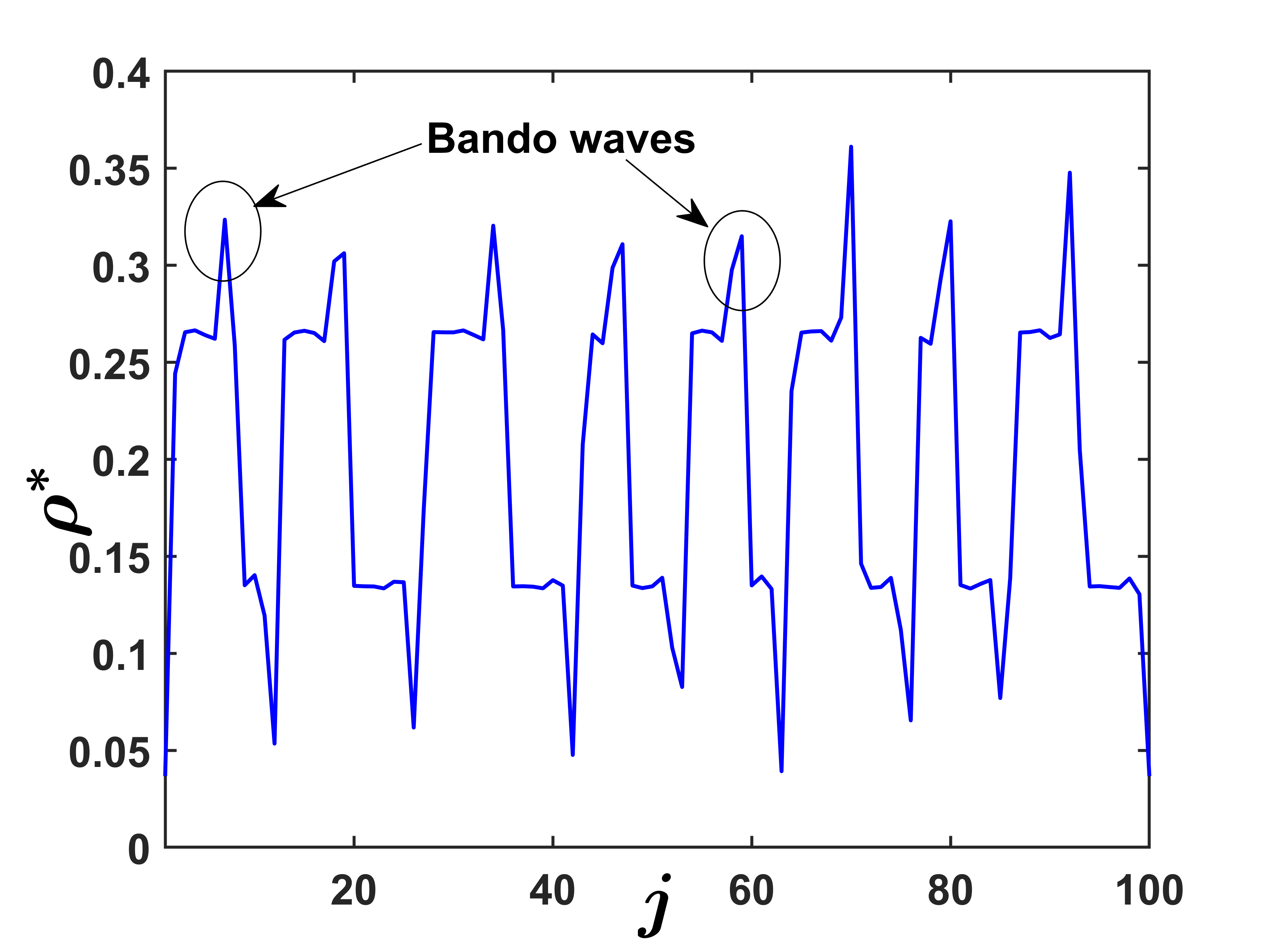} \quad
          \includegraphics[width=5.25cm,height=4.5cm]{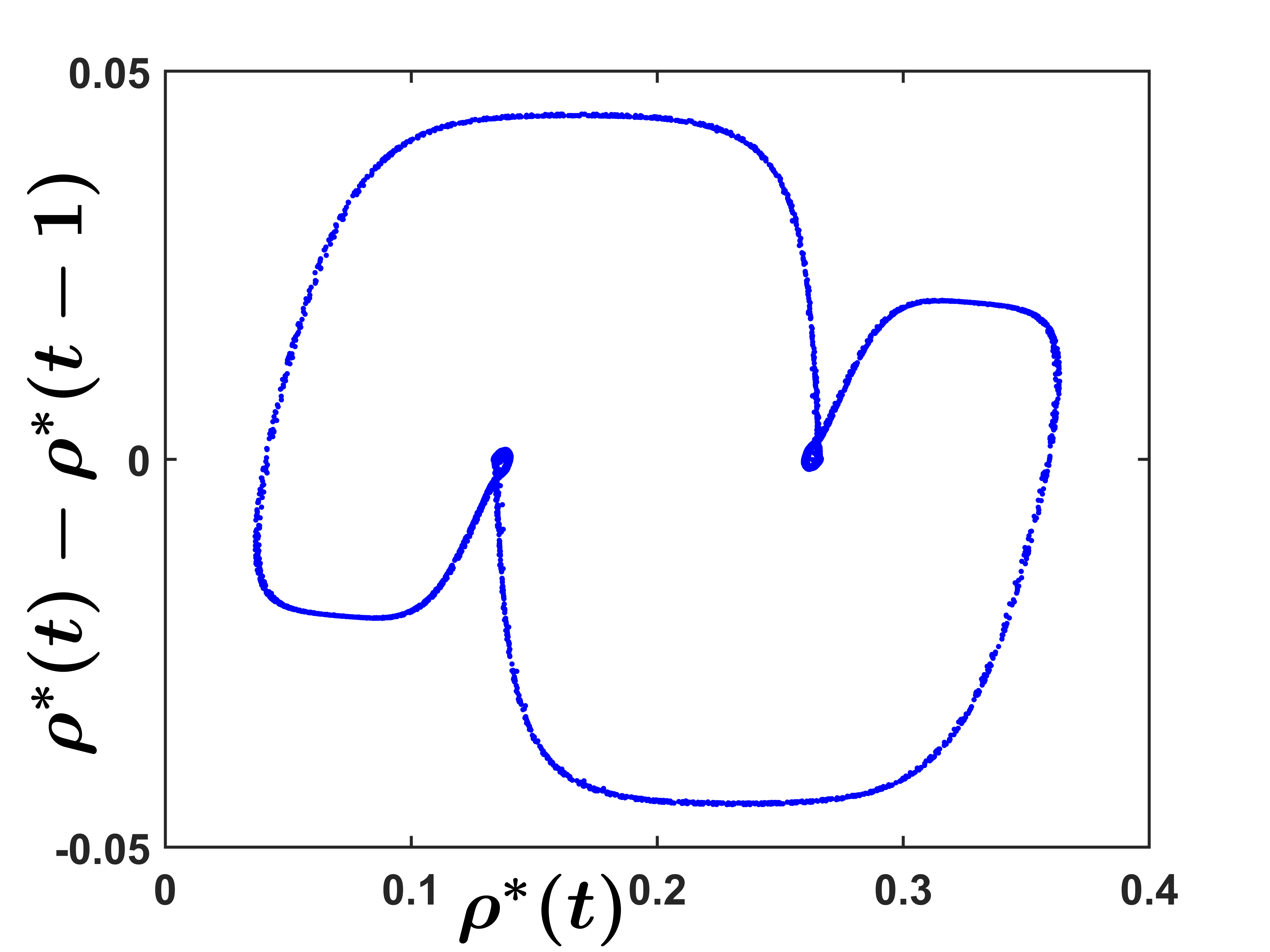} \\
          (a) \hspace{5cm} (b) \hspace{5cm} (c)\\
       \includegraphics[width=5.25cm,height=4.5cm]{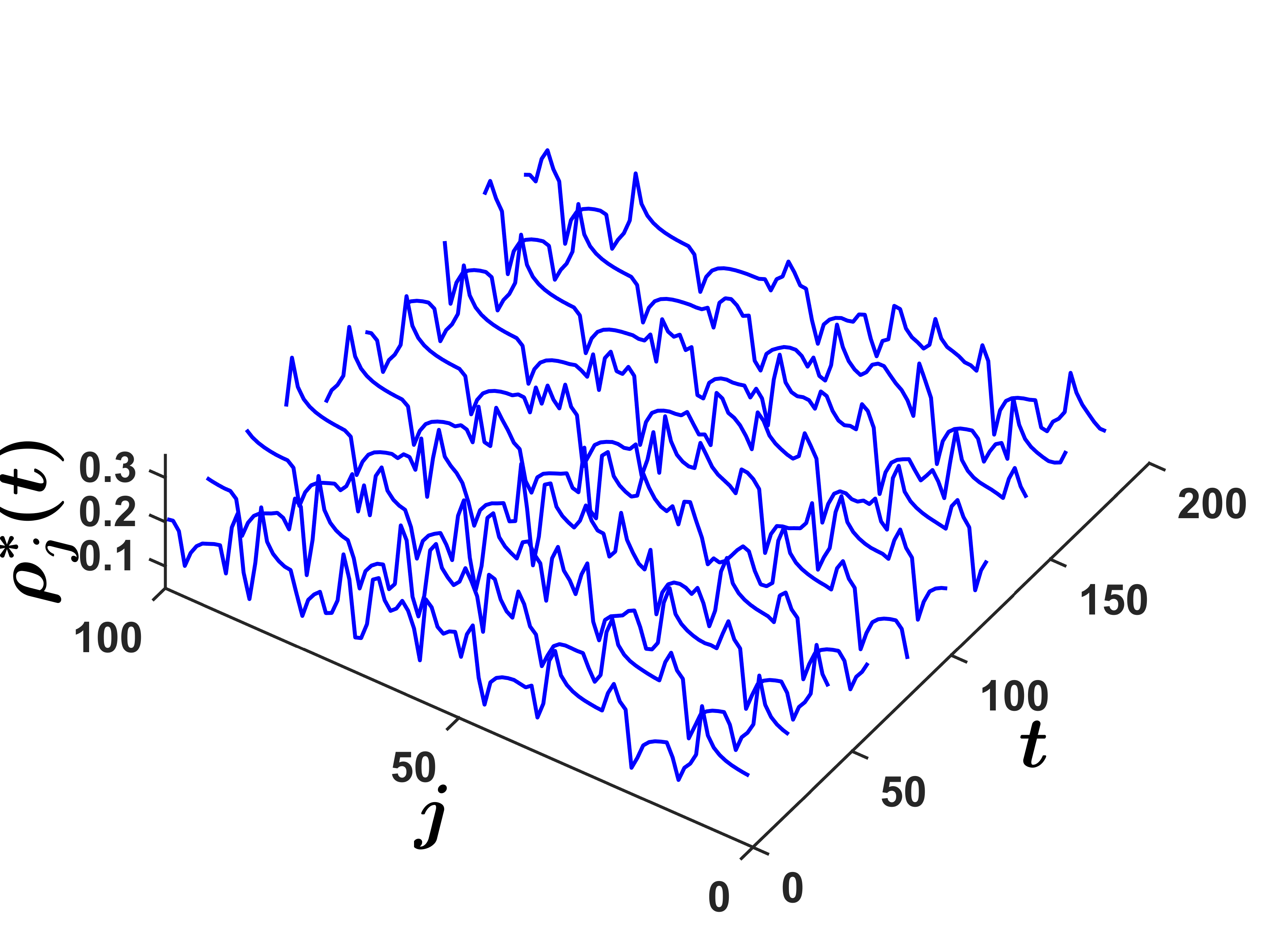}\quad
          \includegraphics[width=5.25cm,height=4.5cm]{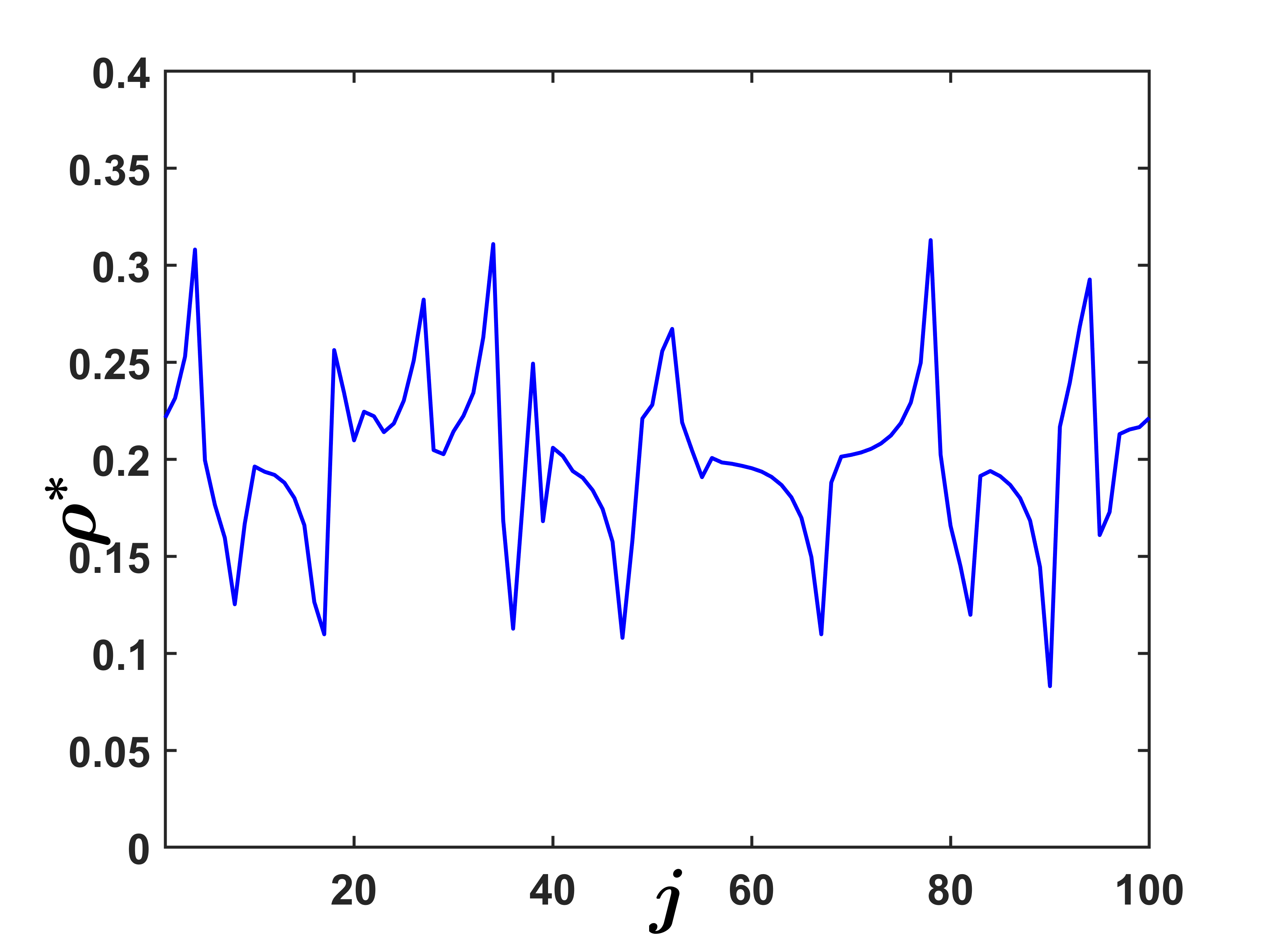}\quad
       \includegraphics[width=5.25cm,height=4.5cm]{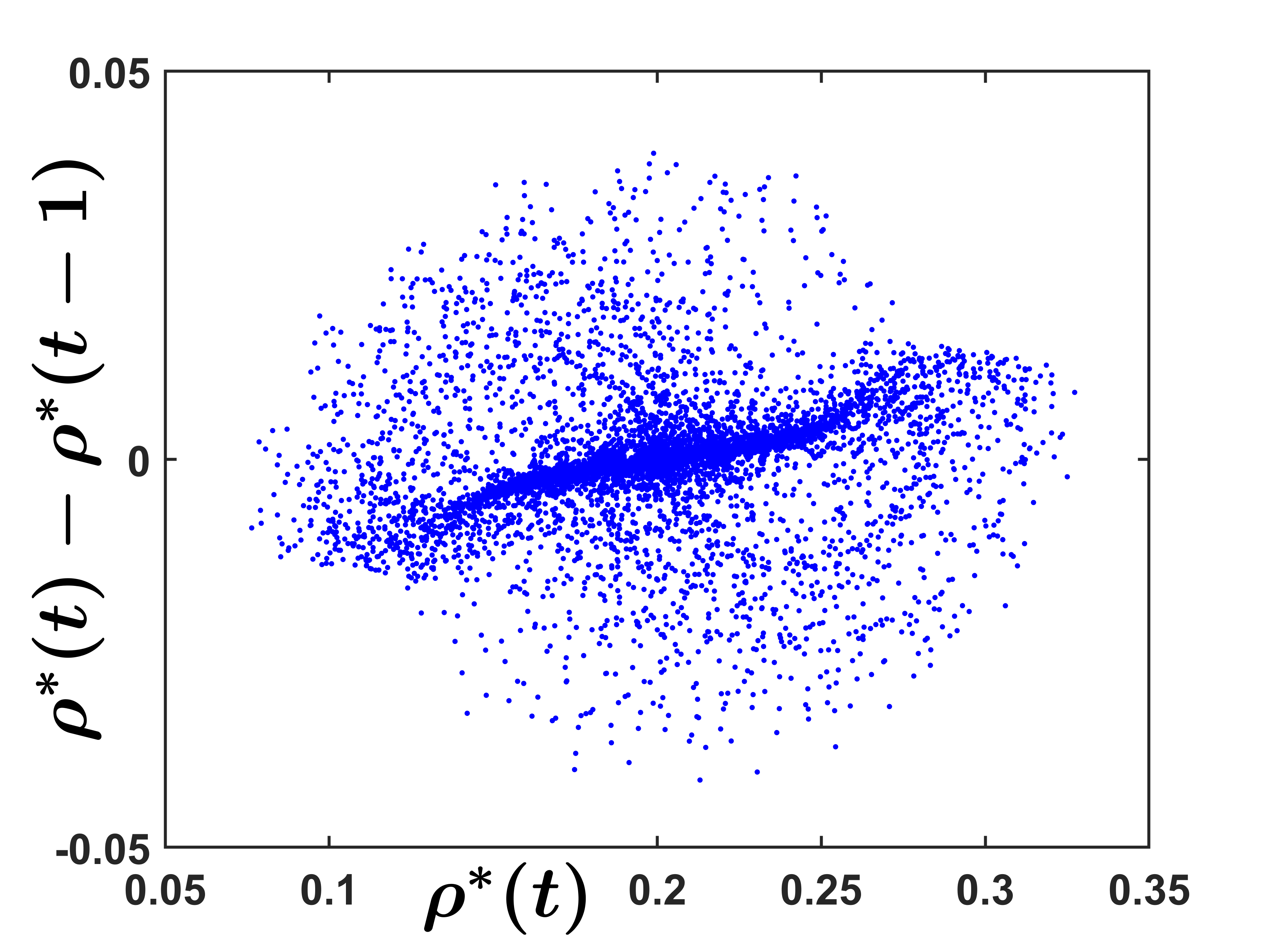}\\ 
       (d) \hspace{5cm} (e) \hspace{5cm} (f)\\
 \caption{Spatiotemporal evolutions of density waves $\rho_j^*(t)$ ($=B\rho_j(t)$) when (a) $a=3.5$, (d) $a=5$, typical density distribution over the $100$  at $N=25200$ simulation time when (b) $a=3.5$, (e) $a=5$, density vs density difference plot between $16000-25200$ simulation time for (c) $a=3.5$, (f) $B=5$. Plots (a), (b), and (c) are associated with kink jam, and (d), (e), and (f) are associated with chaotic jam instability. Other parameters are $B=1.6,$ $C=0.7,$ $\gamma=0.4$ and $\rho_0=\rho_c=0.2.$}
 \label{fig4} 
\end{figure*}
  \begin{figure*}[!t] 
 \centering
 
    \includegraphics[width=16cm,height=10cm]{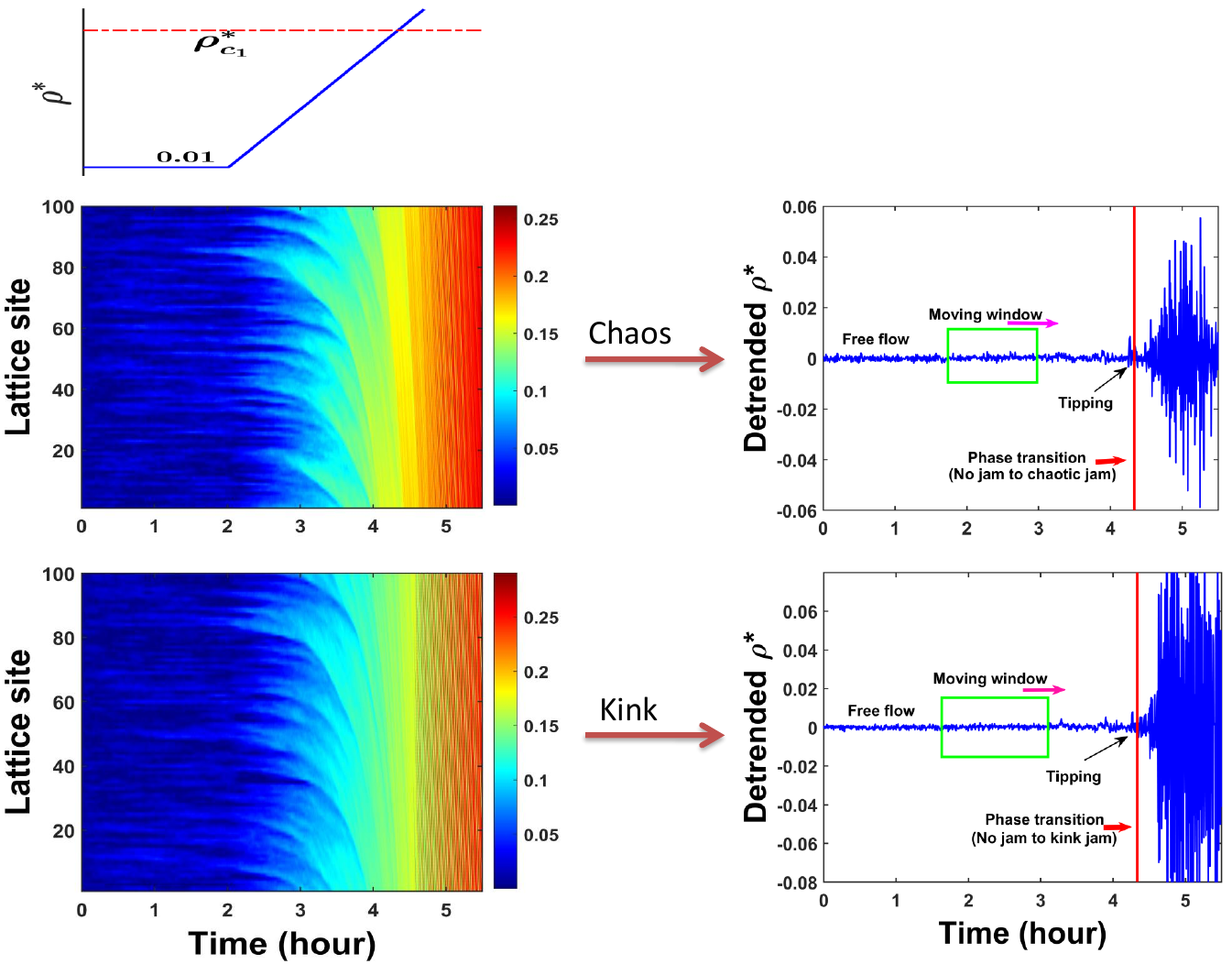}\\
    (a)\\
      \includegraphics[width=16cm,height=4cm]{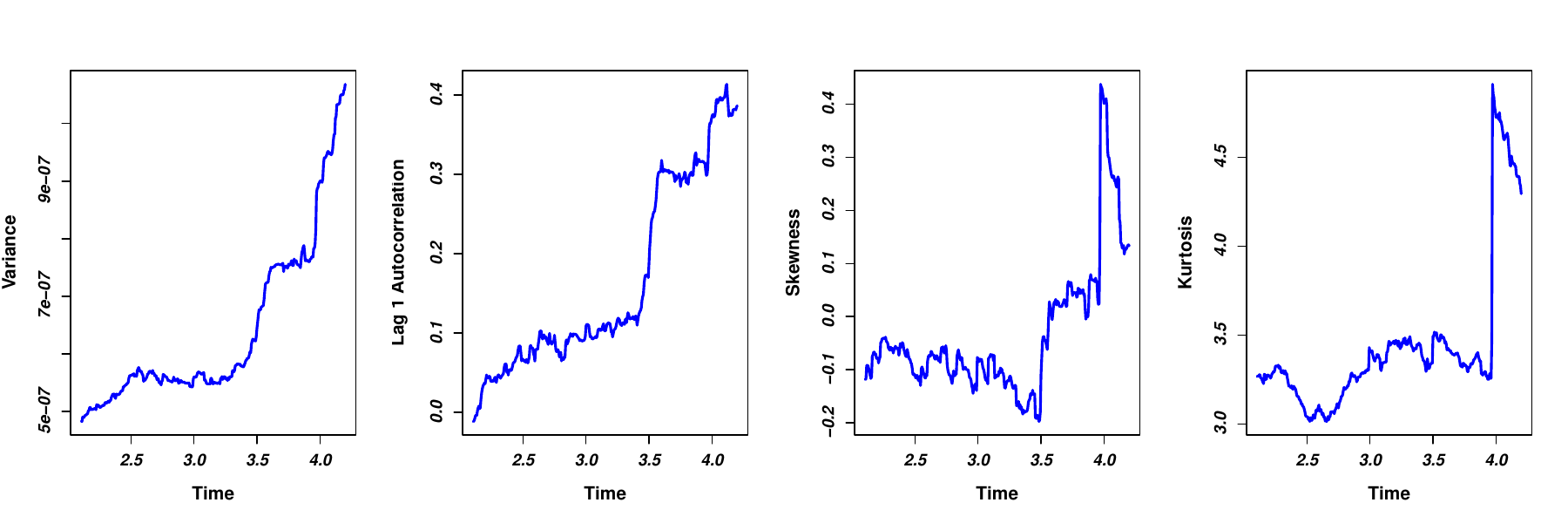}\\
    (b)\\
       \includegraphics[width=16cm,height=4cm]{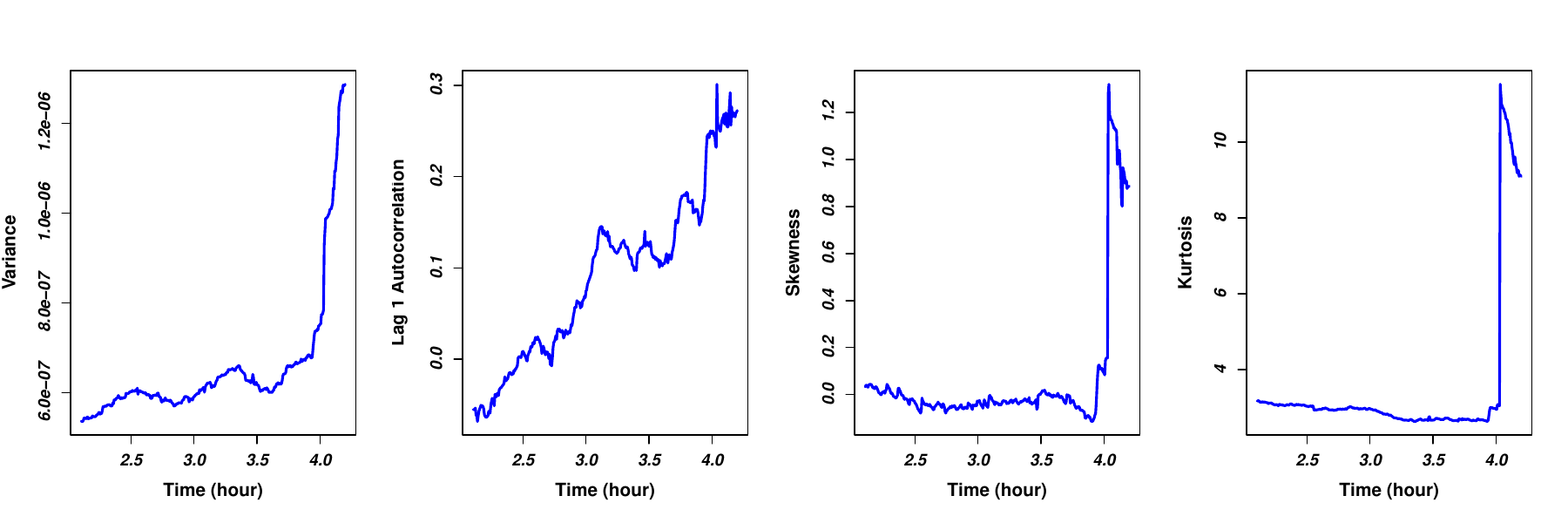}\\
       (c)
        \caption{
        (a) The spatiotemporal pattern of traffic flow density. The mean vehicular density on the road remains constant during the initial two-hour simulation period and afterward exhibits a linear increase until it reaches its critical threshold.  Patterns of chaotic jam and kink jam are detected when the system's dynamics surpasses its stability threshold $\rho_{c_1}^*=0.1573$ for two different values of $a$, namely $a=5$ and $a=3.5$. Detrended density was retrieved from five specific lattice sites located in the center of the system to determine EWSs. The early warning signals were derived from stochastic fluctuations observed in the detrended mean density of vehicles at certain lattice sites, for both the (b) chaotic scenario (c) and the kink scenario.
}
        \label{fig3}
\end{figure*}
\subsection{Early Warning Signals} \label{EWS}
It is apparent that tipping points (TPs) manifest across a wide range of scales since scientific investigations have revealed evidence of TPs occurring at the microbiological level, specifically in yeast \cite{dai2012generic}, as well as on a global scale, leading to climatic shifts \cite{dakos2008slowing}.
The potential future imprints of a forthcoming TP in a dynamical system can be predicted by analyzing the patterns in the statistical moments of time series data, which are referred to as early warning signals (EWSs) \cite{scheffer2009early}.  The stability or instability of a dynamical system is closely linked to the concept of EWSs as they are effective in expressing CSD mathematically. The system's slow rate of returning to a stable state in close proximity to a transition results in significant swings of the system state around the stable state. Consequently, this leads to an increase in statistical measures such as variance and lag-1 autocorrelation. From the statistical point of view variance is the second-order moment of a distribution around the mean $\mu$ and it serves as an early warning indicator measured as standard deviation \cite{carpenter2006rising}. CSD can also be detected by the change in the correlation structure of a time series derived from a dynamical system approaching instability. Lag-1 autocorrelation is the simplest measure to detect the relationship between consecutive observations. As the research area corresponding to TP for different dynamical systems is expanding rapidly, researchers are exploring new methods for detecting abrupt drift of the system from time series data \cite{proverbio2023systematic,deb2022machine,o2023ewsmethods} as well as from experimental data \cite{veraart2012recovery,boers2017deforestation}.
When there is a possibility of a traffic jam a driver senses it beforehand and intuitively decelerates his speed to adjust to the traffic condition. In the vicinity of a traffic jam, the speeds of vehicles slow down systematically and the rate of overall dynamics gradually decreases. As CSD happens, the system under consideration becomes more correlated with its past and thus autocorrelation estimated from the time series gradually increases. Perturbations can accumulate, which gives rise to the size of fluctuations and as a consequence variance or other higher order statistical indicators increase \cite{scheffer2012anticipating,carpenter2006rising}. For our investigation of the disordered traffic system, we have chosen moving window variance, lag-1 autocorrelation, skewness, and kurtosis to anticipate upcoming instability. Statistically, skewness measures asymmetry, and kurtosis measures tailedness of a distribution. We consider $L=100$ lattices with periodic boundary conditions for simulation. The dynamics are subsequently collected at intervals of $20$ seconds in order to obtain statistical information. The LH model for disordered traffic is constituted of two control parameters, namely density and velocity. Based upon the values of any of these control parameters a nearly free traffic flow can become congested. From linear stability analysis, we have obtained two critical densities $\rho^*_{c_1}$ and  $\rho^*_{c_2}$ (see $Appendix$ \ref{Appendix}), in between them free flow loses its linear stability and becomes vulnerable to small perturbations as they grow over time, and lead to congestion. We use generic early warning indicators to predict traffic instability before its occurrence by analyzing whether the statistical metrics (EWSs) are sensitive to impending regime shift or not.

 \section{Results}\label{Results}
\subsection{Numerical Simulation of the Traffic Flow Model} \label{NST}

Now, we carry out extensive numerical simulations to validate our theoretical findings from linear as well as nonlinear stability analysis and to generate time series data for testing the effectiveness of generic EWSs. We graphically illustrate theoretical results to understand the impact of area occupancy factor $B$ and passing rate $\gamma$ in the case of heterogeneous disordered traffic represented by Eq. (\ref{eq:5v2}). Fig. \ref{fig2}(a) represents the neutral boundaries (solid lines) as obtained from linear stability analysis (\ref{LSA}) and the coexisting curves (dotted lines) as obtained from nonlinear stability analysis (\ref{NLSA})  for different values of $B$ with  $\gamma=0.05$ and $C=0.7$. The coexisting curves divide the phase plane into three subregions namely stable, metastable, and unstable regions. From Fig. \ref{fig2}(a) it is clear that an increment of area occupancy factor $B$ uplifts the neutral boundaries as well as the coexisting curves, which means that an increase in disorder can effectively reduce the stable region of the traffic flow. From nonlinear stability analysis (\ref{NLSA}) we obtain the critical value $\gamma<\frac{1}{14}$ for the existence of coexisting curves. Now we investigate the effect of $\gamma$ shown in Fig. \ref{fig2}(b) where $B=1.6$ and $C=0.7.$ Clearly, the coexisting curve exits for $\gamma=0.05$ only, which validates the result obtained from nonlinear stability analysis. It is evident that passing rate $\gamma$ has a positive effect on the expansion of the unstable region. Next, we plot the sensitivity diagram Fig. \ref{fig2}(c) in $(\rho^*, B, a)$ space for different values of $\gamma$ to illustrate the combined effect of disorder and passing on neutral stability boundary in case of unidirectional traffic flow and we observe that area occupancy factor and passing plays a significant role to increase the instability. In Fig. \ref{fig2}(d) the unstable region in the $(\gamma, a)$ phase plane is further divided into kink jam region and chaotic jam region where the line $a_c=-\frac{7B\rho_0^2CV^\prime(B\rho_0)}{2}$ acts as the separatrix between the two different regions. From Fig. \ref{fig2}(d) we conclude that kink jam will occur only if the condition $a < a_c$ is satisfied otherwise chaotic jam will occur.\\To understand the distinct features of the two aforementioned cases, we simulate the density evolution equation for unidirectional disordered traffic allowing passing without on/off ramps represented by Eq. (\ref{eq:5v2}) for $L=100$ lattices with periodic boundary conditions. The simulation is performed for a sufficiently long time to obtain the stationary density distribution. Initially, traffic densities are chosen as uniform $\rho_0=0.2$ except near the central lattice, where a small perturbation $\Delta \rho =0.05$ is given as
  $$ \rho^*_j(0)=\rho^*_j(1)=
 \begin{cases} 
      \rho_0; & j \ne \frac{L}{2}-1,\frac{L}{2} \\
       \rho_0+\Delta \rho; & j = \frac{L}{2}-1 \\
    \rho_0-\Delta \rho; & j = \frac{L}{2}.\\
   \end{cases}
 $$
\textbf{Case 1}: When $a<a_c$:\\
\hspace{5cm} 
Fig. \ref{fig4}(a) shows the space time evolution of traffic density for $B=1.6,$ $\gamma=0.4,$ $C=0.7$ and $\rho_c=0.2$ with sensitivity $a(=\frac{1}{\tau})=3.5$ and $t=25000-25200.$ It is a typical traffic pattern of kink-antikink type as $a<a_c=3.93.$ This kink-antikink traffic waves are regular and robust in the stationary state. Kink density waves propagate in the backward direction and represent the kink jam.
Fig. \ref{fig4}(b) shows typical density distribution over all the $100$ lattices at stationary at stationary state. This density profile has two propagating waves with different speeds, separated by an expanding and contracting zone of congested traffic density. The spikes appearing on the nonlinear traffic waves are commonly referred to as Bando waves \cite{gupta2014analyses}  and are highlighted in circles. Kink jam exhibits periodic motion and to show this, we plot the phase space diagram of density difference $\rho^*(t)-\rho^*(t-1)$ against $\rho^*(t)$ for $t=16000-25200$ in Fig. \ref{fig4}(c). The observed data points exhibit a clustering pattern within a nearly elliptical ring, suggesting that kink jam exhibits regular periodic orbit behavior like a \vspace{0.07cm}limit cycle.\\
\textbf{Case 2}: When $a>a_c$:\\
To illustrate this situation, we simulate the traffic system (\ref{eq:5v2}) while keeping all the parameters same as that of Case 1 except $a=5$ and from the theoretical investigation we are expected to get chaotic density waves as $a>a_c=3.93$. When we plot the spatiotemporal evolution of density waves after a sufficiently large time in Fig. \ref{fig4}(d), we observe the appearance of irregular and aperiodic density waves signifying chaos. Further
Fig. \ref{fig4}(e) shows a typical distribution of these chaotic density waves over all the $100$ lattices at stationary state.  These density waves propagate in the backward direction and represent the chaotic jam. In order to justify the `chaotic' term we plot the phase space diagram of density difference $\rho^*(t)-\rho^*(t-1)$ against $\rho^*(t)$ for $t=16000-25200$ in Fig. \ref{fig4}(f) and it is actually a set of dispersed points in a bounded domain like a strange attractor and therefore chaotic traffic waves exhibit similarities to the concept of chaos as described in dynamical system theory \vspace{0.07cm}\cite{strogatz2018nonlinear}.\\Hence, it can be inferred that the theoretical results align with the outcomes obtained from numerical simulations. The variables denoting disorder in the traffic system and the passing rate exhibit significant influence on the collective traffic behavior, thereby resulting in diverse stationary density profiles.  For higher values of sensitivity, the phenomenon of phase transition is observed from uniform flow $\longrightarrow$ chaotic jam $\longrightarrow$ uniform flow. Conversely, for lower values of sensitivity, the phase transition happens from uniform flow $\longrightarrow$ kink jam $\longrightarrow$ uniform flow.

\subsection{Analysis of EWSs Extracted from Simulated Data} \label{ResultsSD}
A traffic flow can advance towards instability depending upon any of these factors affecting the traffic dynamics, like mean vehicular density, bad weather conditions, specific time of day, sudden accidents or emergency situations, etc. The prediction of traffic instability or congestion is a crucial element that aids in comprehending the intricate nonlinear dynamics of vehicular traffic systems and offers useful insights for developing appropriate control strategies within the field of intelligent transportation systems. In order to identify preliminary indicators of traffic disruptions relevant in emerging economies we have generated time series data from the mathematical model represented by Eq. (\ref{eq:5}). The mean vehicular density on the road remains constant at 0.01 for the initial two hours and afterward, it is linearly increased over time so that the system can reach into the vicinity of critical threshold $\rho^*_{c_1}=0.1573$ (see the left top panel of Fig. \ref{fig3}(a)), which is determined for a certain choice of parameters (see $Appendix$ (\ref{LSA})). Our theoretical findings and numerical simulations suggest that depending upon the chosen parameter values density waves in the traffic jam can be of kink or chaotic type. Fig. \ref{fig3}(a) demonstrates the spatiotemporal evolution and detrended average density profile of  $5$ lattice sites in the middle portion of the road (namely $48$-th, $49$-th, $50$-th, $51$-st and $52$-nd)  for kink jam as well as chaotic jam scenario where the parameters are same as Case 1 and Case 2 of (\ref{NST}) respectively. These detrended density profiles act as required time series data to retrieve the statistical information.\\
For the extraction of EWSs, we monitor the stochastic fluctuations due to Gaussian white noise for each $20$ second in the time series data for both types of instabilities. We introduce variance, lag-1 autocorrelation, skewness, and kurtosis as our guiding metric to predict incoming TP in advance. The second-order moment variance refers to the mathematical expectation of the squared deviations from the mean. The first panel of  Fig. \ref{fig3}(b) and Fig. \ref{fig3}(c) shows the variance for chaotic and kink instability and we observe it increases gradually as the system approaches instability implying the increment of variability along the mean toward instability. Autocorrelation refers to the measurement of the correlation between a signal and a delayed replica of itself, where the correlation is assessed as a function of the time delay. In the present investigation, we examine the lag-1 autocorrelation to understand the correlation structure between values separated by a single time step. The second panel of Fig. \ref{fig3}(b) and Fig. \ref{fig3}(c) demonstrates a gradual increase in the lag-1 autocorrelation metric as the system moves toward chaotic and kink instability respectively. The increment in lag-1 autocorrelation implies an augmentation of the short-term memory capacity of the traffic system because of CSD \cite{pavithran2021effect}.  The last two panels of Fig. \ref{fig3}(b) and Fig. \ref{fig3}(c) indicate the skewness and kurtosis respectively obtained from the available time series data.  Skewness is a statistical measure that quantifies the degree of symmetry in the probability distribution of a dataset relative to its mean while kurtosis provides insights into whether the probability distribution of the data has large tails or is more centrally concentrated in relation to a normal distribution. We observe that variance and lag-1 autocorrelation increase systematically whereas in the vicinity of congestion skewness and kurtosis show a sudden jump.\\
Clearly, all the introduced metrics are showing a concurrent rise in their values and hence they are effective in predicting traffic jams beforehand for a heterogeneous disordered traffic system. It is evident that rising EWSs occur before the actual regime shift between stable free flow to unstable congested flow for kink as well as chaotic instability. It is posited that EWSs do not exhibit sensitivity to the type of instability and thus are effective in anticipating any type of traffic jam in advance which is a matter of practical significance.

%
\section{Discussion and Conclusion} \label{DC}
Different real-world phenomena like epidemiological outbreaks, population collapse, traffic congestion, etc can switch to an alternate stable state in case of multistability or can become unstable in case of monostability due to crossing of tipping point. Particularly in the case of traffic flow, a nearly free flow can become congested when it crosses the critical density threshold due to a sudden change in a parameter value that affects the overall traffic dynamics. In this study, we proposed an LH area occupancy
model with passing to examine the dynamics of a heterogeneous disordered traffic system, which is particularly significant in the context of emerging economies where road infrastructure often lacks lane discipline. We observed that heterogeneity and passing play a crucial role in shaping the overall traffic dynamics.\\
Linear stability analysis yields two critical values for density, such that traffic flow is unstable if $\rho_{c_1}<\rho < \rho_{c_2}$. From nonlinear stability analysis, we further split the unstable region into two parts namely kink and chaotic region. We considered a closed hypothetical road for simulation.  We took the initial density far less from the lower critical threshold $\rho_{c_1}$ and after some initial phase of stationary simulations, density is poured linearly through an on-ramp on the closed loop to anticipate both kink and chaotic jam. Noise is ubiquitous in physical or biological systems and it plays a significant role in shaping the dynamics \cite{sharma2015stochasticity}. Gaussian white noise in our system plays the role of external perturbations. We extracted data from a segment consisting of five lattices and found that noise accumulated as the system approached instability. As a direct consequence, the statistical indicators are showing a sharp hike, warning us about an impending regime shift. The evaluations conducted in the text concentrated on the lower critical density, which is a prevalent kind of instability observed on roadways. However, it is possible to achieve comparable outcomes by assessing the suggested approaches for the critical transition at the higher critical density. \\ The potential future scope of the proposed research involves investigating the resilience of early warning signals and extending their applicability to traffic networks.\\
In conclusion, the early warning signals enable the anticipation of traffic congestion in advance within a heterogeneous and disordered unidirectional traffic flow. The potential outcomes of our study can be effectively applied by gathering data from a heterogeneous disordered highway and tuning the parameters accordingly.

\section{Appendix} \label{Appendix}
\subsection{Linear Stability Analysis} \label{LSA}
The purpose of doing linear stability analysis is to investigate the response of traffic flow to perturbations in the flow system when it is under a steady state. The solution of the uniform steady state is $\rho_j(t)=\rho_0$ and $v_j(t)=V(\rho_0).$ Let $y_j(t)$ represent a minor deviation from the uniform steady-state flow, which leads to $\rho_j(t)=\rho_0+y_j(t).$ Substituting this perturbed density profile into Eq. (\ref{eq:5v2}), we get the following linear equation
\begin{equation}\label{eq:16}
    \Delta y_j(t+\tau) +\tau \rho_0^2 BCV^\prime(B\rho_0)\bigl(\Delta y_j(t)-\gamma \Delta^2 y_j(t)\bigr)=0,
    \end{equation}
    
where 
  $$\begin{cases} 
      \Delta_{t} y_{j}(t+\tau)=y_{j}(t+2\tau)-y_{j}(t+\tau),\\
       \Delta y_j(t)=y_{j+1}(t)-y_j(t),\\
     \Delta^{2} y_j(t)=y_{j+2}(t)-2y_{j+1}(t)+y_j(t),\\ 
   \end{cases}$$
 and $V^\prime(B\rho_0)=\frac{dV(B\rho_j)}{d\rho_j}\biggr|_{\rho_j=\rho_0}$. 
 By expanding $y_j(t)=e^{(ikj+wt)}$, the following equation of $w$ is derived from Eq. (\ref{eq:16})
  \begin{align}
     \begin{split}
  e^{2w\tau}-e^{w\tau}&+\tau \rho_0^2 BC\bigl(e^{ik}-1\bigr)V^\prime(B\rho_0)\\&-\gamma \tau \rho_0^2BC\bigl(e^{ik}-1\bigr)^2V^\prime(B\rho_0) =0.
        \end{split}
    \label{eq:17}
 \end{align}
 Inserting $w=w_1(ik)+w_2(ik)^2+...$ into Eq. (\ref{eq:17}), the first and second order terms that are the coefficients $ik$ and $(ik)^2$ are obtained as follows
  \begin{equation}  \label{eq:18}
  \begin{split}
   w_1 &= \displaystyle{-\rho_0^2BCV^\prime(B\rho_0),} \\[5pt]
   w_2 &=-\frac{3}{2}\tau w_1^2-\frac{(1-2\gamma)}{2}\rho_0^2BCV^\prime(B\rho_0).\\
   \end{split}
 \end{equation}
  If $w_2$ is a negative value, the uniform steady state flow becomes unstable. When $w_2$ is a positive value, the uniform flow is stable. The neutral stability condition is thus achieved when $w_2=0.$ Thus,
   neutral stability curve is as follows in terms of sensitivity $a(=\frac{1}{\tau})$  
 \begin{equation}
   a=-\frac{(1-2\gamma)}{3\rho_0^2BCV^\prime(B\rho_0)}
    \label{eq:19}
 \end{equation}
 For small disturbances of long wavelengths, the uniform traffic flow is stable if
  \begin{equation}
   a > -\frac{3\rho_0^2BCV^\prime(B\rho_0)}{(1-2\gamma)}.\\
    \label{eq:21}
 \end{equation}
 If we further write $\rho^*=B\rho_0,$ then 
the stability condition (\ref{eq:21}) is associated with two critical values of density, denoted as $\rho^*_{c_1}$ and $\rho^*_{c_2}$, with $\rho_{c_1}^*<\rho_{c_2}^*$. 
For $B=1.6$, $C=0.7$, $\gamma=0.4$ and $a=3.93$ the values of $\rho_{c_1}^*$ and $\rho_{c_2}^*$ are given by $\rho_{c_1}^*=0.1573$ and $\rho_{c_2}^*=0.2743$.
\begin{table*}[!t]
 \centering
  \caption{Values of $h_i$ coefficients}
\begin{tabular}{ |c| c| c|c|c|  } 
 \hline
  $h_1$ & $h_2$ & $h_3$ & $ h_4$ &$h_5$ \\ 
  \hline
  
  $b+BC\rho_0^2V^\prime$& $\frac{3b^2\tau}{2}+\frac{(1-2\gamma)}{2}BC\rho_0^2V^\prime$ & $  \frac{B^2C\rho_0^2V^{\prime\prime}}{2}$ &$ \frac{7b^3\tau^2}{6} +\frac{BC\rho_0^2V^{\prime}}{6}-\gamma BC\rho_0^2V^{\prime}$&$ \frac{B^2C\rho_0^2V^{\prime\prime}(1-2\gamma)}{4}$
 \\ \hline

   $ h_6$ &   $ h_7$ & $ h_8$ &  $ h_9$ & $ h_{10}$ \\ \hline 
   $\frac{B^3C\rho_0^2V^{\prime \prime \prime}}{6}$&   $3b\tau$ & $\frac{5b^4\tau^3}{8}+\frac{BC\rho_0^2V^{\prime}}{24}-\frac{7\gamma BC\rho_0^2V^{\prime}}{12}$ &
       $ \frac{B^2C\rho_0^2V^{\prime\prime}(1-6\gamma)}{12}$ & $ \frac{B^3C\rho_0^2V^{\prime\prime\prime}(1-2\gamma)}{12}$ \\ \hline
\end{tabular}
    \label{Table2}
\end{table*}
  \begin{table*}[!t]
 \centering
 \caption{Values of $b_i$ coefficients}
\begin{tabular}{ |c| c| c|c|c|  } 
 \hline
  $b_1$ & $b_2$ & $b_3$& $ b_4$ & $ b_5$\\
  \hline
  $-\frac{7}{6}b^3\tau_c^2-\frac{1}{6}BC\rho_c^2V^\prime+BC\rho_c^2V^\prime\gamma$ & $\frac{B^3C\rho_c^2V^{\prime\prime\prime}}{6}$ & $\frac{3b^2\tau_c}{2}$ &
  $\frac{5}{8}b^4\tau_c^3+\frac{BC\rho_c^2V^{\prime}}{24}-\frac{7}{12}BC\rho_c^2V^{\prime}\gamma+3b\tau_cb_1$& $\frac{B^3C\rho_c^2V^{\prime\prime\prime}(1-2\gamma)}{12}-3b\tau_cb_2$\\
  \hline
\end{tabular}
    \label{Table1}
\end{table*}
\subsection{Nonlinear Stability Analysis} \label{NLSA}
In order to get insight into the gradual changes occurring in the vicinity of the critical point $(\rho_c,a_c)$, we introduce two variables $X$ and $T$ to represent the slow dynamics. These variables are associated with a minute positive scaling parameter $\epsilon$, where $\epsilon$ is significantly smaller than $1$. 
\begin{equation}
    \left.
     \begin{array}{ccccc}
  \displaystyle{ X }&=& \displaystyle{\epsilon(j+bt),} \\[5pt] 
  \displaystyle{T} &=&\displaystyle{\epsilon^3t,}\\
        \end{array}
  \right.
    \label{eq:23}
 \end{equation}
 where the unknown parameter $b$ is needed to be found.
The assumption about density is 
\begin{equation}
    \left.
     \begin{array}{ccccc}
  \displaystyle{ \rho_j(t) }&=& \displaystyle{\rho_c+\epsilon R(X,T).} 
        \end{array}
  \right.
    \label{eq:13}
 \end{equation}
 From Eqs. (\ref{eq:5v2}), (\ref{eq:23}), (\ref{eq:13}), using Taylor series expansion to the fifth order of $\epsilon$, we finally obtain
 \begin{equation}
  \epsilon^2h_1\partial_XR +P\epsilon^3 +Q\epsilon^4+R\epsilon^5= 0.
    \label{eq:32} 
 \end{equation}
 where  $$\begin{cases} 
      P= \bigl(h_2\partial_X^2R +h_3\partial_XR^2 \bigr) ,\\
    Q= \bigl(\partial_TR +h_4\partial_X^3R+h_5\partial_x^2R^2+h_6\partial_XR^3\bigr),\\
     R=\bigl(h_7\partial_T\partial_XR+h_8\partial_X^4R +h_9\partial_X^3R^2+h_{10}\partial_X^2R^3 \bigr).
   \end{cases}$$
 The coefficients of the Eq. (\ref{eq:32}) is given in Table (\ref{Table2}).
 By substituting $ b=-BC\rho_0^2V^\prime(\rho_c^*)$  where $\rho_c^*=B\rho_c$ and $\tau=\tau_c(1+\epsilon^2)$ into Eq. (\ref{eq:32}), we rewrite it as
 \begin{align}
     \begin{split}
  \epsilon^4(\partial_TR-b_1\partial_X^3R+b_2\partial_XR^3)+\epsilon^5(&b_3\partial_X^2R+b_4\partial_X^4R\\&+b_5\partial_X^2R^3)=0.
        \end{split}
        \label{eq:25}
 \end{align}

Using this relation $T=\frac{T^{\prime}}{b_1}$ and $R=\sqrt{\frac{b_1}{b_2}}R^{\prime}$ we transform the Eq. (\ref{eq:25}) into the regularized equation given below
\begin{align}
     \begin{split}
  \partial_T^{\prime}R^{\prime}-\partial_X^3R^{\prime}+\partial_XR^{\prime 3}+\frac{\epsilon}{b_1}\biggl(b_3\partial_X^2R^{\prime}&+b_4\partial_X^4R^{\prime}\\&\hspace{-1cm}+\frac{b_1 b_5}{b_2}\partial_X^2R^{\prime 3}\biggr)=0. 
        \end{split}
    \label{eq:26}
 \end{align}
 Our transformation is thus given by
  \begin{equation}
    \left.
     \begin{array}{ccccc}
  \displaystyle{ T^\prime}&=& \displaystyle{\biggl(-\frac{(1-13\gamma-14\gamma^2)(BC\rho_c^2V^\prime)}{27}\biggr)T,} \\[5pt]
    \displaystyle{  R}&=& \displaystyle{\left(-\frac{2(1-13\gamma-14\gamma^2)(BC\rho_c^2V^\prime)}{9B^3C\rho_c^2V^{\prime \prime \prime}}\right)^{\frac{1}{2}}R^\prime,} 
        \end{array}
  \right.
    \label{eq:30}
 \end{equation}
 where we assumed that 
 \begin{equation}
    \left.
     \begin{array}{ccccc}
  \displaystyle{ 13\gamma+14\gamma^2}&<& \displaystyle{1.} 
        \end{array}
  \right.
    \label{eq:31}
 \end{equation}
 The coefficients of the Eq. (\ref{eq:25}) is given by the Table (\ref{Table1}).
 Ignoring the \begin{huge}o\end{huge}$(\epsilon)$ term in Eq. (\ref{eq:26}), we get mKdv equation whose desired kink soliton solution is given as follows
  \begin{equation}
 R_0^{\prime}(X,T^{\prime})=\sqrt{\mu}\left[tanh\left(\sqrt{\frac{\mu}{2}}(X-\mu T^{\prime})\right)\right],
    \label{eq:27}
 \end{equation}
 where, $\mu$ is the propagation velocity and can be determined by solving the condition: \\$(R_0,M[R_0])=\int_{-\infty}^{\infty}dXR_0M[R_0]=0,$ where $$M[R_0]=\frac{1}{b_1}\left(b_3\partial_X^2R^{\prime}+b_4\partial_X^4R^{\prime}+\frac{b_1 b_5}{b_2}\partial_X^2R^{\prime 3}\right).$$\\
  The suitable choice of $\mu$ is obtained as 
  \begin{equation}
     \begin{array}{ccccc}
  \displaystyle{  \mu}&=& \displaystyle{\frac{5b_2b_3}{2b_2b_4-3b_1b_5}}.
        \end{array}
    \label{eq:28}
 \end{equation}
  Hence the kink solution is 
 \begin{equation}
       \begin{array}{ccccc}
  \displaystyle{ \rho_j^*}&=& \displaystyle{\rho_c+\alpha\biggl[tanh\left(\sqrt{\frac{\mu}{2}}(X-\mu b_1T)\right)\biggr]} ,
        \end{array}
    \label{eq:29}
 \end{equation}
 with the amplitude $\alpha=\sqrt{\frac{b_1\epsilon^2 \mu}{b_2}}$ and $\epsilon^2=\left(\frac{\tau}{\tau_c}-1\right)$.
\\ The kink solution exists only if condition (\ref{eq:31}) is satisfied. So the kink solution exists if 
 \begin{equation}
     0 \leq \gamma  < \frac{1}{14}
     \label{eq:75}
 \end{equation}
 If the value of $\gamma$ exceeds the threshold $\frac{1}{14}$, then the modified KdV equation cannot be derived from the nonlinear analysis. If $\gamma >\frac{1}{14}$ then the scaling assumption which we have employed in Eq. (\ref{eq:23}) breaks down. So to describe the chaotic density waves we need a different scaling assumption.

%

\vfill
\end{document}